\documentclass[journal, final, onecolumn, 10pt]{IEEEtran}
\usepackage[compress]{cite}

\usepackage{subfigure}
\usepackage{float}

\usepackage[bookmarks=false, hidelinks]{hyperref}

\usepackage{amsmath,amssymb,amsfonts}
\allowdisplaybreaks
\usepackage{bm}                         
\usepackage{multirow}
\usepackage{cases}                       
\usepackage{tikz}

\usepackage[normalem]{ulem}             

\usepackage[linesnumbered,noend,ruled,vlined]{algorithm2e}
\usepackage{algpseudocode}

\usepackage{graphicx}
\usepackage{textcomp}
\usepackage{xcolor}

\usepackage{balance}

\begin{document}

\title{Extreme Distance Distributions \\ of Poisson Voronoi Cells
\thanks{
}
}


\author{\IEEEauthorblockN{Jaume Anguera Peris, Joakim Jaldén}\\ 
\IEEEauthorblockA{School of Electrical Engineering and Computer Science\\
KTH Royal Institute of Technology, Stockholm, Sweden\\
Email: \{jaumeap, jalden\}@kth.se}
}

\maketitle
\begin{abstract}
Poisson point processes provide a versatile framework for modeling the distributions of random points in space. When the space is partitioned into cells, each associated with a single generating point from the Poisson process, there appears a geometric structure known as Poisson Voronoi tessellation. These tessellations find applications in various fields such as biology, material science, and communications, where the statistical properties of the Voronoi cells reveal patterns and structures that hold key insights into the underlying processes generating the observed phenomena.

In this paper, we investigate a distance measure of Poisson Voronoi tessellations that is emerging in the literature, yet for which its statistical and geometrical properties remain explored only in the asymptotic case when the density of seed points approaches infinity. Our work, specifically focused on homogeneous Poisson point processes, characterizes the cumulative distribution functions governing the smallest and largest distances between the points generating the Voronoi regions and their respective vertices for an arbitrary density of points in $\mathbb{R}^2$. For that, we conduct a Monte-Carlo type simulation with $10^8$ Voronoi cells and fit the resulting empirical cumulative distribution functions to the Generalized Gamma, Gamma, Log-normal, Rayleigh, and Weibull distributions. Our analysis compares these fits in terms of root mean-squared error and maximum absolute variation, revealing the Generalized Gamma distribution as the best-fit model for characterizing these distances in homogeneous Poisson Voronoi tessellations. Furthermore, we provide estimates for the maximum likelihood and the $95$\% confidence interval of the parameters of the Generalized Gamma distribution along with the algorithm implemented to calculate the maximum and minimum distances.

\end{abstract}


\section{Introduction}
\label{sec:introduction}

\IEEEPARstart{F}{rom} the distribution of stars in galaxies and the arrangement of atoms in crystal lattices to the intricate realms of cellular biology, communication networks, and material science, we often encounter the occurrence of spatial points representing the locations of various phenomena or objects of interest. The arrangement of these points in space not only reflects the underlying processes generating them but also reveals patterns and structures that hold key insights into the nature of the phenomena being studied. Specially, when the study of these spatial points extends beyond simple observation and involves the analysis of structures or parameters subject to randomness, stochastic geometry provides a mathematical framework to analyze the spatial patterns of such random objects. Within this framework, the concept of spatial point processes emerges, providing a formalized model for characterizing the random distribution of points in space.

In the realm of spatial point processes, we encounter Poisson point processes (PPP), characterized by their spatial randomness and independence, making them powerful tools for modeling a diverse array of phenomena. One of the most studied spatial structures derived from PPPs is the Voronoi tessellation. Voronoi tessellations partition the space into regions based on the proximity of points in space to the so-called generating points or seeds. Specifically, each point in space is assigned to the nearest generating point, resulting in a partitioning of the space into cells, where each cells contains one generating point as well as all the points that are closer to that generating point than to any other. Within the resulting geometric structure, we can distinguish $n+1$ types of points in $\mathbb{R}^n$ \cite{brakke1987statistics}. For the two-dimensional case, a point closest to only one generating point belongs to the interior of a cell. A point equidistant to two generating points lays on the boundary between two regions, and a point equidistant to three generating points corresponds to the vertex where three cells intersect.

These simple, yet powerful concepts, have led to the significant theoretical analysis of the distance between generator points and vertices \cite[Eq.~6]{muche2005poisson}, the distance between generator points \cite[Eq.~19]{muche2005poisson}, or the angles between neighbouring generators points and vertices \cite[Eq.~21]{muche2005poisson}, to name a few. However, due to the stochastic nature of Poisson point processes, deriving closed-form analytical expressions for certain mathematical properties of Voronoi tessellations can be challenging or even intractable. These challenges, coupled with high relevance of such analytical expressions, often requires the use of numerical methods, simulations, and mathematical approximations to characterize the statistical properties of Poisson Voronoi tessellations. Examples of these are the Voronoi's cell area, perimeter, or number of vertices \cite{ferenc2007size,tanemura2003statistical}.

Altogether, thanks to all these prior works, Poisson Voronoi tessellations have been used to cluster biological data~\cite{edla2012clustering}, analyze the territorial land between animals~\cite{stewart2010voronoi}, characterize the spread of infections~\cite{bherwani2021understanding}, identify clusters of stars and galaxies~\cite{ramella2001finding}, classify binding pockets in proteins~\cite{feinstein2021bionoi}, or model millimeter wave cellular systems~\cite{andrews2016modeling}. Nevertheless, among the unknown geometrical and statistical properties of Poisson Voronoi tessellations, the maximum and minimum distances between the generator points and the vertices of their Vornoi cells are distance measures that have only been explored before in the asymptotic case when the density of generator points in space approaches infinity \cite{calka2014extreme}. Yet, investigating the statistical properties of these distance measures under finite point densities holds significant practical relevance in modeling the worse-case user/node in mobile communications to ensure reliable signal coverage while minimizing interference and signal degradation~\cite{peris2022modelling};  analyzing energy densities of periodic ferromagnetic Ising systems to minimize anisotropic surface tensions~\cite{chambolle2023crystallinity}; 
analyzing interatomic distances to optimize the synthesis processes, mechanical properties, and performance of various materials \cite{jalem2018general}; characterizing the signal transduction pathways and biomolecular interaction of cells, organelles, and biomolecules \cite{xia2016review}; or optimizing the spatial placement of emergency shelters to increase safety in urban environments \cite{ma2019emergency}.

\subsection{Our contribution}
In this paper, we specifically focus on homogeneous PPP, which correspond to a specific instance of PPPs where the intensity or rate of point occurrence is constant across the entire spatial domain. Within homogeneous PPP, we further refine our focus to characterizing the statistical properties of the normalized maximum and minimum distances between a generator seed and the vertices of its Voronoi cell. Since we are interested in finite distances, we examine only those cells that are entirely contained within the tessellation and exclude outer cells. Considering this, we conduct a Monte-Carlo type simulation and approximate the resulting cumulative distribution functions (CDFs) of the maximum and minimum distances to the Generalized Gamma, Gamma, Log-normal, Rayleigh, and Weibull distributions, similar to \cite{ferenc2007size,tanemura2003statistical,gezer2021statistical}. We compare the fit distribution to the empirical distributions in terms of the root mean-squared error and maximum absolute variation. We show that the Generalized Gamma distribution is the best-fit for characterizing the maximum and minimum distances in homogeneous Poisson Voronoi tessellations and we derive the maximum likelihood and 95\% confidence interval of its parameters.



\subsection{Document organization}
The remainder of this paper is organized as follows. Section \ref{sec:ProblemFormulation} formulates the problem at hand and defines the quantities of interest in a graphical manner. In Section \ref{sec:ExtremeDistDistributions}, we divide the analysis of the maximum and minimum distances in two parts. First, we calculate the theoretical CDFs of the maximum and minimum distance in $1$D to gain insights on the relationship between the parameters of interest and the stochastic processes governing the Poisson Voronoi diagrams. Then, we discuss the mathematical foundations to derive the CDF of the maximum and minimum distances in 2D and illustrate the best-fit approximation alongside the empirical CDF obtained from the Monte-Carlo simulation. We conclude the paper in Section \ref{sec:Conclusions}, followed by an Appendix with the comparison between fit functions and the insights on the algorithms we used for calculating the maximum and minimum distances.

\subsection{Notation}
Throughout this paper, sets are represented by calligraphic letters $\mathcal{S}$, random variables are represented by upper-case letters $X$, and samples (observations) of random variables are represented by lower-case letters $x$. Specially, we reserve $D$ and $R$ to refer to the random variables representing Euclidean distances in one-dimensional and two-dimensional spaces, respectively. The cumulative distribution function of a continuous random variable $X$ is represented by $F_X(x) = \mathbb{P}(X\leq x)$.
For the random samples $x_1, x_2, \dots, x_N$ drawn independently for a continuous distribution $F_X(x)$, we denote by $x_{(i)}$ the $i$-th order statistic of the sample, defined by sorting the samples $x_1, x_2, \dots, x_N$ in increasing order. For an integer $N$, $[N]$ stands for $\{1,2,\dots, N\}$.

\begin{figure*}
    \centering
    \subfigure[]{\includegraphics[width = 0.47\textwidth]{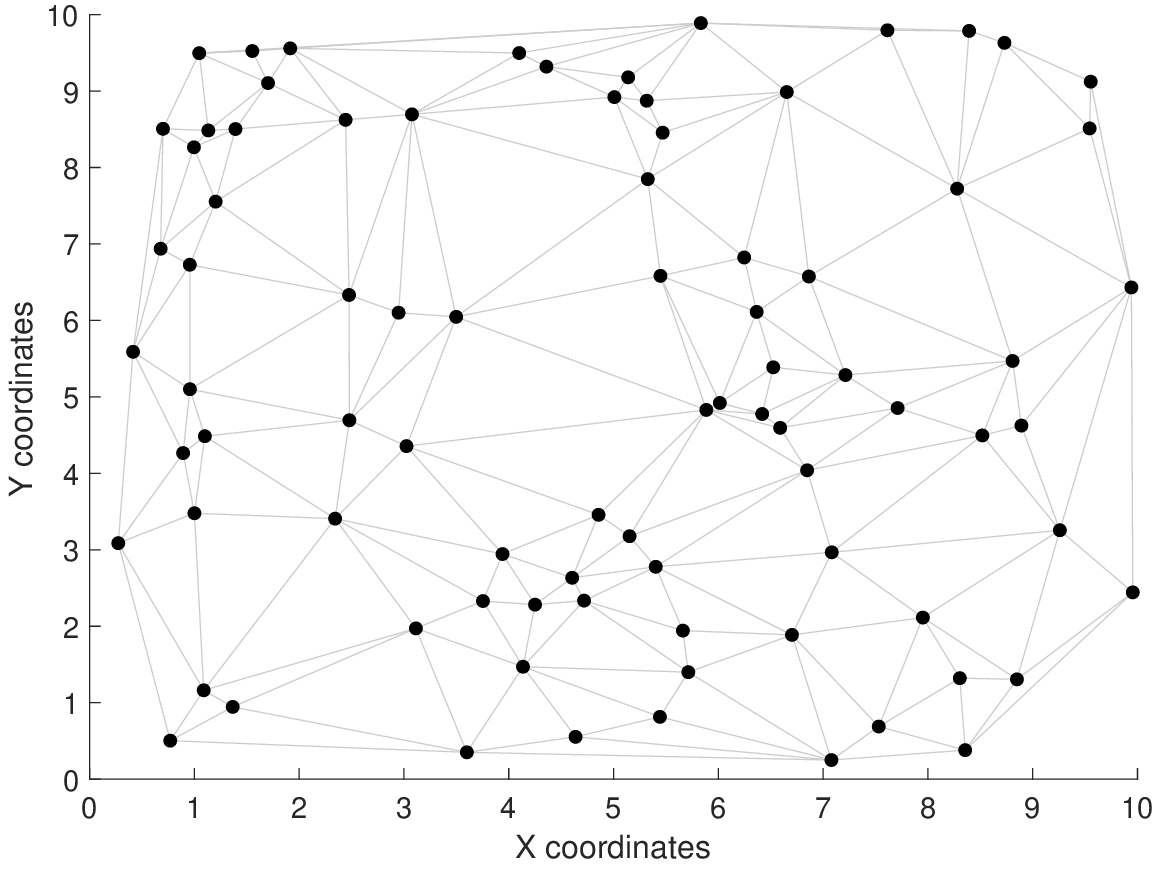}}$\quad$
    \subfigure[]{\includegraphics[width = 0.47\textwidth]{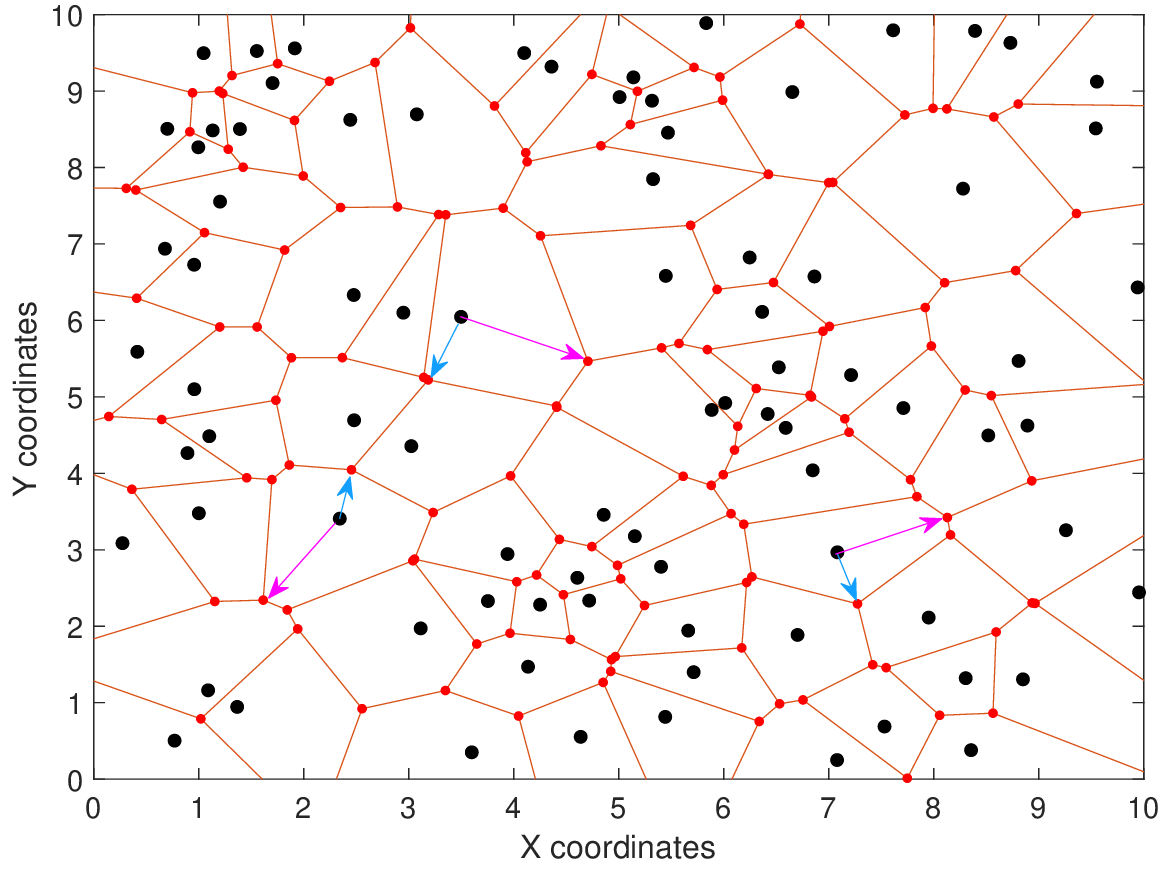} \label{fig:VoronoiTessellation_withArrows} }
    \caption{(a) Delanuay triangulation of a sample realization of an homogeneous PPP with intensity parameter $\lambda=0.8$ over an area $A=100$. (b) Poisson Voronoi tessellation of the homogeneous PPP generated in (a). Black dots represent generator seeds in the planar space. Red dots represent the vertices of the Voronoi cells, which correspond to the circumcenters of the triangles generated by the Delanuay triangulation. Gray lines delimit the triangles of the Delanuay triangulation, and orange lines represent the edges of the Voronoi cells. The arrows in blue and magenta exemplify the stochastic quantities of interest in $2$D.}
    \label{fig:VoronoiTessellation_2D}
\end{figure*}

\section{Problem formulation}
\label{sec:ProblemFormulation}
Consider an homogeneous PPP in $\mathbb{R}^2$ with intensity parameter $\lambda > 0$, and denote the countably infinite set of \emph{generator  seeds} in the spatial domain as $\mathcal{P} = \{p_1, p_2, \dots, p_N\}$. At the same time, envision a Voronoi tessellation generated from the Delanuay triangulation of the planar point set $\mathcal{P}$, as shown in Figure \ref{fig:VoronoiTessellation_2D}. Within this tessellation, define the Voronoi region of a seed point $p_i$, denoted by $\mathcal{V}(p_i)$, as
\begin{equation*}
    \mathcal{V}(p_i) = \Big\{ x\in \mathbb{R}^2 : \|p_i - x \|_2 \leq \|p_j - x \|_2, \; \forall p_j\in\mathcal{P}\setminus p_i \Big\},\quad \forall p_i\in\mathcal{P}.
\end{equation*}
Notice that the Voronoi tessellation is a geometric construction that partitions the space into regions based on the proximity of points in the space to the seed points $p_i\in\mathcal{P}$, and the Voronoi diagram of $\mathcal{P}$, defined as the union of all Voronoi regions, $\mathcal{V} = \cup_{i=1}^N \mathcal{V}(p_i)$, constitutes the entire space domain.

For any two generator seeds, $p_i$ and $p_j$, define the \emph{edge} between their associated Voronoi regions $\mathcal{V}(p_i)$, $\mathcal{V}(p_j)$, as
\begin{equation*}
    \mathcal{A}(p_i,p_j) = \mathcal{V}(p_i) \cap \mathcal{V}(p_j), \quad \forall j\neq i.
\end{equation*}
If an edge satisfies $\mathcal{A}(p_i,p_j)\neq \emptyset$, the Voronoi regions $\mathcal{V}(p_i)$ and $\mathcal{V}(p_j)$ are considered to be adjacent. All adjacent cells to a Voronoi region are considered to be neighbouring cells. Considering this, let the set of indices of seed points of neighbouring cells to the Voronoi region $\mathcal{V}(p_i)$ be defined as
\begin{equation*}
    \mathcal{N}(p_i) = \Big\{ j\in[N] \setminus i : \mathcal{A}(p_i,p_j) \neq \emptyset \Big\}.
\end{equation*}

Furthermore, refer as \emph{vertex} the points in the plane that are equidistant to three generator seeds. From the point of view of a Voronoi region $\mathcal{V}(p_i)$, a vertex is constituted by the endpoints of the edges of two neighbouring cells that are also adjacent to each other. Mathematically speaking, the vertices of a Voronoi region $\mathcal{V}(p_i)$ are defined as
\begin{equation*}
    \ell(p_i) = \Big\{ x\in\mathbb{R}^2 : \|p_i - x \|_2 = \|p_j - x \|_2 = \|p_k - x \|_2, \, \exists\, j,k \in \mathcal{N}(p_i), \, j\neq k \neq i \Big\}.
\end{equation*}

Considering the above, we are interested in obtaining relevant insights about the maximum and minimum distances in the Voronoi diagram. In essence, the maximum and minimum distances of a Voronoi region $\mathcal{V}(p_i)$ correspond to the smallest and largest Euclidean distances between the generator seed $p_i$ and any of its vertices, i.e.,
\begin{align*}
    r_\mathrm{min}(p_i) = \min_{x\in\ell(p_i)} \| x - p_i\|_2, \quad \text{and} \quad r_\mathrm{max}(p_i) = \max_{x\in\ell(p_i)} \| x - p_i\|_2,\quad \forall p_i\in \mathcal{P}.
\end{align*}
Notice, however, that in the process of creating a Voronoi tessellation, the initial step involves the stochastic placement of the generator seed points within the spatial domain. As a consequence of this randomness in the spatial arrangement of generator seeds, the Voronoi tessellation and the properties of Voronoi cells results to be stochastic in nature. Therefore, in order to characterize the maximum and minimum distances in the entire Voronoi diagram, we must account for the diverse realizations of the Voronoi tessellation, and consider that individual Voronoi regions lead to varying maximum and minimum distances. In this regard, the next section focuses on characterising the maximum and minimum distances in Poisson Voronoi diagrams by means of their statistical distributions.


\section{Maximum and minimum distance distributions}
\label{sec:ExtremeDistDistributions}
To better understand the problem at hand, let us first study theoretically the CDF of the maximum and minimum distances in the one-dimensional case, and let us also present a graphical example to gain insights on the relationship between the distribution functions and the characteristics of the underlying homogeneous PPP generating the seed points.

\subsection{1-dimensional case}
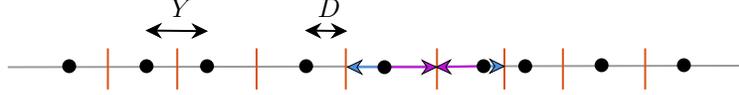
\begin{figure*}
\centering
    \tikzset{every picture/.style={line width=0.75pt}} 

\begin{tikzpicture}[x=0.75pt,y=0.75pt,yscale=-1,xscale=1]

\draw [color={rgb, 255:red, 155; green, 155; blue, 155 }  ,draw opacity=0.79 ]   (100,109.5) -- (473,108.75) ;
\draw  [fill={rgb, 255:red, 0; green, 0; blue, 0 }  ,fill opacity=1 ] (128,109) .. controls (128,107.34) and (129.34,106) .. (131,106) .. controls (132.66,106) and (134,107.34) .. (134,109) .. controls (134,110.66) and (132.66,112) .. (131,112) .. controls (129.34,112) and (128,110.66) .. (128,109) -- cycle ;
\draw  [fill={rgb, 255:red, 0; green, 0; blue, 0 }  ,fill opacity=1 ] (247.5,109) .. controls (247.5,107.34) and (248.84,106) .. (250.5,106) .. controls (252.16,106) and (253.5,107.34) .. (253.5,109) .. controls (253.5,110.66) and (252.16,112) .. (250.5,112) .. controls (248.84,112) and (247.5,110.66) .. (247.5,109) -- cycle ;
\draw  [fill={rgb, 255:red, 0; green, 0; blue, 0 }  ,fill opacity=1 ] (197.5,109) .. controls (197.5,107.34) and (198.84,106) .. (200.5,106) .. controls (202.16,106) and (203.5,107.34) .. (203.5,109) .. controls (203.5,110.66) and (202.16,112) .. (200.5,112) .. controls (198.84,112) and (197.5,110.66) .. (197.5,109) -- cycle ;
\draw  [fill={rgb, 255:red, 0; green, 0; blue, 0 }  ,fill opacity=1 ] (167,109) .. controls (167,107.34) and (168.34,106) .. (170,106) .. controls (171.66,106) and (173,107.34) .. (173,109) .. controls (173,110.66) and (171.66,112) .. (170,112) .. controls (168.34,112) and (167,110.66) .. (167,109) -- cycle ;
\draw  [fill={rgb, 255:red, 0; green, 0; blue, 0 }  ,fill opacity=1 ] (396.5,108.5) .. controls (396.5,106.84) and (397.84,105.5) .. (399.5,105.5) .. controls (401.16,105.5) and (402.5,106.84) .. (402.5,108.5) .. controls (402.5,110.16) and (401.16,111.5) .. (399.5,111.5) .. controls (397.84,111.5) and (396.5,110.16) .. (396.5,108.5) -- cycle ;
\draw  [fill={rgb, 255:red, 0; green, 0; blue, 0 }  ,fill opacity=1 ] (358,109) .. controls (358,107.34) and (359.34,106) .. (361,106) .. controls (362.66,106) and (364,107.34) .. (364,109) .. controls (364,110.66) and (362.66,112) .. (361,112) .. controls (359.34,112) and (358,110.66) .. (358,109) -- cycle ;
\draw  [fill={rgb, 255:red, 0; green, 0; blue, 0 }  ,fill opacity=1 ] (438,108.5) .. controls (438,106.84) and (439.34,105.5) .. (441,105.5) .. controls (442.66,105.5) and (444,106.84) .. (444,108.5) .. controls (444,110.16) and (442.66,111.5) .. (441,111.5) .. controls (439.34,111.5) and (438,110.16) .. (438,108.5) -- cycle ;
\draw [color={rgb, 255:red, 217; green, 83; blue, 25 }  ,draw opacity=1 ]   (150.5,100) -- (150.5,120.75) ;
\draw [color={rgb, 255:red, 217; green, 83; blue, 25 }  ,draw opacity=1 ]   (185.5,100) -- (185.5,120.75) ;
\draw [color={rgb, 255:red, 217; green, 83; blue, 25 }  ,draw opacity=1 ]   (225.5,100) -- (225.5,120.75) ;
\draw [color={rgb, 255:red, 217; green, 83; blue, 25 }  ,draw opacity=1 ]   (270.5,100) -- (270.5,120.75) ;
\draw [color={rgb, 255:red, 217; green, 83; blue, 25 }  ,draw opacity=1 ]   (316.5,100) -- (316.5,120.75) ;
\draw [color={rgb, 255:red, 217; green, 83; blue, 25 }  ,draw opacity=1 ]   (350.5,100) -- (350.5,120.75) ;
\draw [color={rgb, 255:red, 217; green, 83; blue, 25 }  ,draw opacity=1 ]   (380,100) -- (380,120.75) ;
\draw [color={rgb, 255:red, 217; green, 83; blue, 25 }  ,draw opacity=1 ]   (421.5,100) -- (421.5,120.75) ;
\draw [line width=0.75]    (173,91.03) -- (197.43,91.29) ;
\draw [shift={(200.43,91.32)}, rotate = 180.61] [fill={rgb, 255:red, 0; green, 0; blue, 0 }  ][line width=0.08]  [draw opacity=0] (7.14,-3.43) -- (0,0) -- (7.14,3.43) -- (4.74,0) -- cycle    ;
\draw [shift={(170,91)}, rotate = 0.61] [fill={rgb, 255:red, 0; green, 0; blue, 0 }  ][line width=0.08]  [draw opacity=0] (7.14,-3.43) -- (0,0) -- (7.14,3.43) -- (4.74,0) -- cycle    ;
\draw [line width=0.75]    (253.43,91.04) -- (267.71,91.24) ;
\draw [shift={(270.71,91.29)}, rotate = 180.81] [fill={rgb, 255:red, 0; green, 0; blue, 0 }  ][line width=0.08]  [draw opacity=0] (7.14,-3.43) -- (0,0) -- (7.14,3.43) -- (4.74,0) -- cycle    ;
\draw [shift={(250.43,91)}, rotate = 0.81] [fill={rgb, 255:red, 0; green, 0; blue, 0 }  ][line width=0.08]  [draw opacity=0] (7.14,-3.43) -- (0,0) -- (7.14,3.43) -- (4.74,0) -- cycle    ;
\draw [color={rgb, 255:red, 189; green, 16; blue, 224 }  ,draw opacity=1 ][line width=0.75]    (291,109.5) -- (313.2,109.5) ;
\draw [shift={(316.2,109.5)}, rotate = 180] [fill={rgb, 255:red, 189; green, 16; blue, 224 }  ,fill opacity=1 ][line width=0.08]  [draw opacity=0] (7.14,-3.43) -- (0,0) -- (7.14,3.43) -- (4.74,0) -- cycle    ;
\draw [color={rgb, 255:red, 74; green, 144; blue, 226 }  ,draw opacity=1 ][line width=0.75]    (291,109.5) -- (274.4,109.5) ;
\draw [shift={(271.4,109.5)}, rotate = 360] [fill={rgb, 255:red, 74; green, 144; blue, 226 }  ,fill opacity=1 ][line width=0.08]  [draw opacity=0] (7.14,-3.43) -- (0,0) -- (7.14,3.43) -- (4.74,0) -- cycle    ;
\draw  [fill={rgb, 255:red, 0; green, 0; blue, 0 }  ,fill opacity=1 ] (287,109.5) .. controls (287,107.84) and (288.34,106.5) .. (290,106.5) .. controls (291.66,106.5) and (293,107.84) .. (293,109.5) .. controls (293,111.16) and (291.66,112.5) .. (290,112.5) .. controls (288.34,112.5) and (287,111.16) .. (287,109.5) -- cycle ;
\draw [color={rgb, 255:red, 74; green, 144; blue, 226 }  ,draw opacity=1 ][line width=0.75]    (340.5,109) -- (347.4,109.07) ;
\draw [shift={(350.4,109.1)}, rotate = 180.58] [fill={rgb, 255:red, 74; green, 144; blue, 226 }  ,fill opacity=1 ][line width=0.08]  [draw opacity=0] (7.14,-3.43) -- (0,0) -- (7.14,3.43) -- (4.74,0) -- cycle    ;
\draw [color={rgb, 255:red, 189; green, 16; blue, 224 }  ,draw opacity=1 ][line width=0.75]    (339.5,109) -- (319.2,109.44) ;
\draw [shift={(316.2,109.5)}, rotate = 358.77] [fill={rgb, 255:red, 189; green, 16; blue, 224 }  ,fill opacity=1 ][line width=0.08]  [draw opacity=0] (7.14,-3.43) -- (0,0) -- (7.14,3.43) -- (4.74,0) -- cycle    ;
\draw  [fill={rgb, 255:red, 0; green, 0; blue, 0 }  ,fill opacity=1 ] (337,109) .. controls (337,107.34) and (338.34,106) .. (340,106) .. controls (341.66,106) and (343,107.34) .. (343,109) .. controls (343,110.66) and (341.66,112) .. (340,112) .. controls (338.34,112) and (337,110.66) .. (337,109) -- cycle ;

\draw (190.5,79.38) node   [align=left] {\begin{minipage}[lt]{12.24pt}\setlength\topsep{0pt}
$\displaystyle Y$
\end{minipage}};
\draw (264.5,79.37) node   [align=left] {\begin{minipage}[lt]{12.24pt}\setlength\topsep{0pt}
$\displaystyle D$
\end{minipage}};

\end{tikzpicture}
    \caption{Poisson Voronoi tessellation of a one-dimensional homogeneous PPP, where $\theta$ represents the distance between generator seed points, and $D$ represents the distances between generator seed points and their edges. The arrows in blue and magenta exemplify the stochastic quantities of interest in $1$D.}
    \label{fig:VoronoiTessellation_1D}
\end{figure*}

Begin by considering a line of length $L$, as illustrated in Figure \ref{fig:VoronoiTessellation_1D}, consisting of $N$ generator seed points uniformly and independently distributed along the line. As we contemplate the limiting case when $N\rightarrow\infty$, $L\rightarrow\infty$, and the density of points per unit length $\lambda = N/L$ remains finite, the PPP naturally emerges as a mathematical model to  capture the stochastic nature of point occurrences along the line. In this limiting case, if we focus on the specific case where the number of points in a segment of the line is linearly proportional to both $\lambda$ and the length of the segment, the inherent randomness in the arrangement of seeds emerges as an homogeneous PPP. 

For any given sample realization of an homogeneous PPP with intensity $\lambda$, the Voronoi diagram of the generating seed points divides the line into multiple bounded Voronoi regions and two unbounded Voronoi regions, associated to the seed points located in the interior and at the ends of the line, respectively. For the sake of the contribution of this paper, we disregard the unbounded regions and focus solely on the interior seed points. In this case, let $Y$ be the random variable that represents the Euclidean distance between two neighbouring seed points  when $N\rightarrow\infty$, $L\rightarrow\infty$, and $\lambda = N/L$ remains finite. Since the arrangement of seeds is modelled as an homogeneous PPP, the CDF of $Y$ can be derived from the null probability of a Poisson distribution as follows:
\begin{equation*}
    F_Y(y) = 1 - \mathbb{P} \big( \nexists \, j\in \mathcal{N}(p_i) : \|p_i - p_j\|_2 \leq y, \, \forall p_i \in \mathcal{P} \big) = 1 - e^{-\lambda y}, \quad y \geq 0.
\end{equation*}
Hence, the random variable $Y$ is Exponentially distributed $Y\sim \text{Exp}(\lambda)$. Notice also that the edge between two neighbouring cells is located at the center of the segment connecting their seed points. Considering this, we can further conclude that that random variable that represents the normalized Euclidean distance between a generator seed point and any of its edges, denoted by $\Bar{D} = \lambda Y/2$, is also Exponentially distributed with probability distribution
\begin{equation}
    F_{\Bar{D}}(d) = 1 - e^{- 2 d}, \quad d \geq 0.
    \label{eq:CDF_dist2Edge}
\end{equation}

From there, we can further observe in Figure \ref{fig:VoronoiTessellation_1D} that all bounded Voronoi regions has two neighbouring cells, one to the left, and another one to the right. Since the process is homogeneous and the locations of seed points are independently distributed, the distances from a seed point to its left and right neighbors are statistically identical. That is, there is no inherent directionality or bias in the placement of points. Consequently, the CDFs of the maximum and minimum distances can be derived from the order statistics of the CDF presented in \eqref{eq:CDF_dist2Edge}. In particular, and without loss of generality, let the normalized distances from a seed point to its left and right neighbours be denoted as $\Bar{D}_1$ and $\Bar{D}_2$, respectively, both with probability distribution~\eqref{eq:CDF_dist2Edge}. Moreover, let $\Bar{D}_{(i)}$ be the $i$-th order statistics of $\Bar{D}_1$ and $\Bar{D}_2$, defined by sorting the values of $\Bar{D}_1$ and $\Bar{D}_2$ in increasing order. Then, it follows that the random variables representing the maximum and minimum distances, $\Bar{D}_{\rm{max}}$ and  $\Bar{D}_{\rm{min}}$, have probability distributions
\begin{align}
    F_{\Bar{D}_{\rm{max}}}(d) = \mathbb{P}\left(\Bar{D}_{(2)}\leq d\right) = \mathbb{P}(\max(\Bar{D}_1,\Bar{D}_2)\leq d) & = \big[F_{\Bar{D}}(d)\big]^2 = \Big[1 - e^{- 2 d}\Big]^2, \quad d \geq 0, \label{eq:ExtremeDist_1D_max}\\[.7em]
    F_{\Bar{D}_{\rm{min}}}(d) = \mathbb{P}\left(\Bar{D}_{(1)}\leq d\right) = \mathbb{P}(\min(\Bar{D}_1,\Bar{D}_2)\leq d) & = 1 - \big[1-F_{\Bar{D}}(d)\big]^2 = 1 - e^{- 4 d}, \quad d \geq 0. \label{eq:ExtremeDist_1D_min}
\end{align}


In summary, by considering an infinite line populated with a finite but dense arrangement of generator points, we established the emergence of homogeneous PPPs as an appropriate model to capture the stochastic nature of point occurrences. From there, it was essential to understand how Voronoi cells are constructed from these seed points. That is, cells are intervals along the line, with edges at the midpoint between neighboring seed points. We then leveraged these geometric properties to characterize the random nature of the normalized maximum and minimum distances in the one-dimensional case.

\subsection{2-dimensional case}
Now that we have gained an insight into the statistical properties of the maximum and minimum distances in one-dimensional, bounded Poisson Voronoi tessellations, we are ready to extend our analysis to the two-dimensional case. For that, notice that the cornerstone for deriving the statistical properties of these distances lies in the geometric construct to partition the space into Voronoi regions. For this reason, we shift our focus to the Delaunay triangulation, as it is the algorithm we consider in this paper for generating Voronoi tessellations.

One of the key properties of the Delaunay triangulation is its efficiency in partitioning the set of generator seed points into triangles while satisfying the Delaunay condition. This condition asserts that a triangulation of the generator seed points is considered a Delaunay triangulation if and only if the circumcircle of each triangle in the triangulation does not contain any other points from the given set $\mathcal{P}$ in its interior. This implies that the set of all vertices in the Voronoi tessellation constitutes the circumcenters of every triangle in the triangulation. Considering these properties of the vertices and the Delanuay triangulation, one can present a configuration space and a configuration measure based on the set of generator seeds $\mathcal{P}$ to derive the statistical properties of several quantities in Poisson Voronoi tessellations \cite{brakke1987statistics}. Specially, if $\Bar{R} = \sqrt{\lambda} R$ represents the random variable that characterizes the normalized Euclidean distance between a seed point and any of its vertices, it follows from the analysis in~\cite{brakke1987statistics} that $\Bar{R}$ has Gamma distribution with probability distribution
\begin{equation}
    F_{\Bar{R}}(r) = 1 - \big(1 + \pi r^2 \big)e^{-\pi r^2}.
    \label{eq:DistributionDistance}
\end{equation}

With that, we have characterized the random nature of the spatial distribution of vertices with respect to the generator seeds. However, the stochastic arrangement of seed points also introduces a variability in the number of vertices associated with each seed point. Moreover, as highlighted in \cite{gezer2021statistical}, even though the interior cells in the Voronoi diagram tend to have more vertices than exterior cells, the Voronoi regions with the highest number of vertices are spread throughout the entire space and do not concentrate at any particular location. Therefore, to effectively capture and analyze the distribution of the maximum and minimum distances within such stochastic Voronoi regions, we require of a more nuanced approach beyond direct reliance on order statistics.

Despite acknowledging the above, any theoretical attempts to finding the analytical expressions for the statistical distributions of the maximum and minimum distances in 2D failed, so we instead resorted to conduct a Monte-Carlo simulation and fit the statistical distributions of the maximum and minimum distances to different functions. In particular, we first performed a Monte-Carlo simulation to obtain the empirical CDF of the normalized maximum and minimum distances following Algorithm \ref{alg:VoronoiTesselation_2D} for a total of $10^8$ Voronoi regions, and then approximated the resulting CDFs to the Generalized Gamma, Gamma, Log-normal, Rayleigh, and Weibull distributions. The results of our findings are summarized in Table \ref{tab:GammaFitParams} in the Appendix.

The best-fit approximation in terms of lowest root mean-squared error and lowest maximum absolute variation resulted to be the Generalized Gamma distribution, defined by
\begin{equation}
    F_{\Bar{X}}(x\,;\,a,b,c) = \frac{a\,b^{c/a}}{\Gamma(c/a)}\, \int_0^x t^{c-1}\,e^{-bt^a}dt, \quad x\geq0, \quad a,\,b,\,c > 0,
    \label{eq:GeneralizedGammaDistribution}
\end{equation}
where $\Gamma(\alpha) = \int_0^{\infty} t^{\alpha-1}e^{-t}dt$ is the Gamma function. Notice that \eqref{eq:GeneralizedGammaDistribution} characterizes the statistical distribution of the normalized distances, but our findings can easily be extended to any arbitrary density of seed points $\lambda$ by considering that
\begin{equation*}
    F_{X}(x\,;\,a,b,c) = F_{\Bar{X}}\left(x\sqrt{\lambda}\,;\,a,b,c\right), \quad \forall \lambda > 0,
\end{equation*}
where $X$ represents the maximum and minimum distances, and $\Bar{X}$ represents the normalized maximum and minimum distances.

\begin{figure*}
    \centering
    \includegraphics[width = 0.47\textwidth]{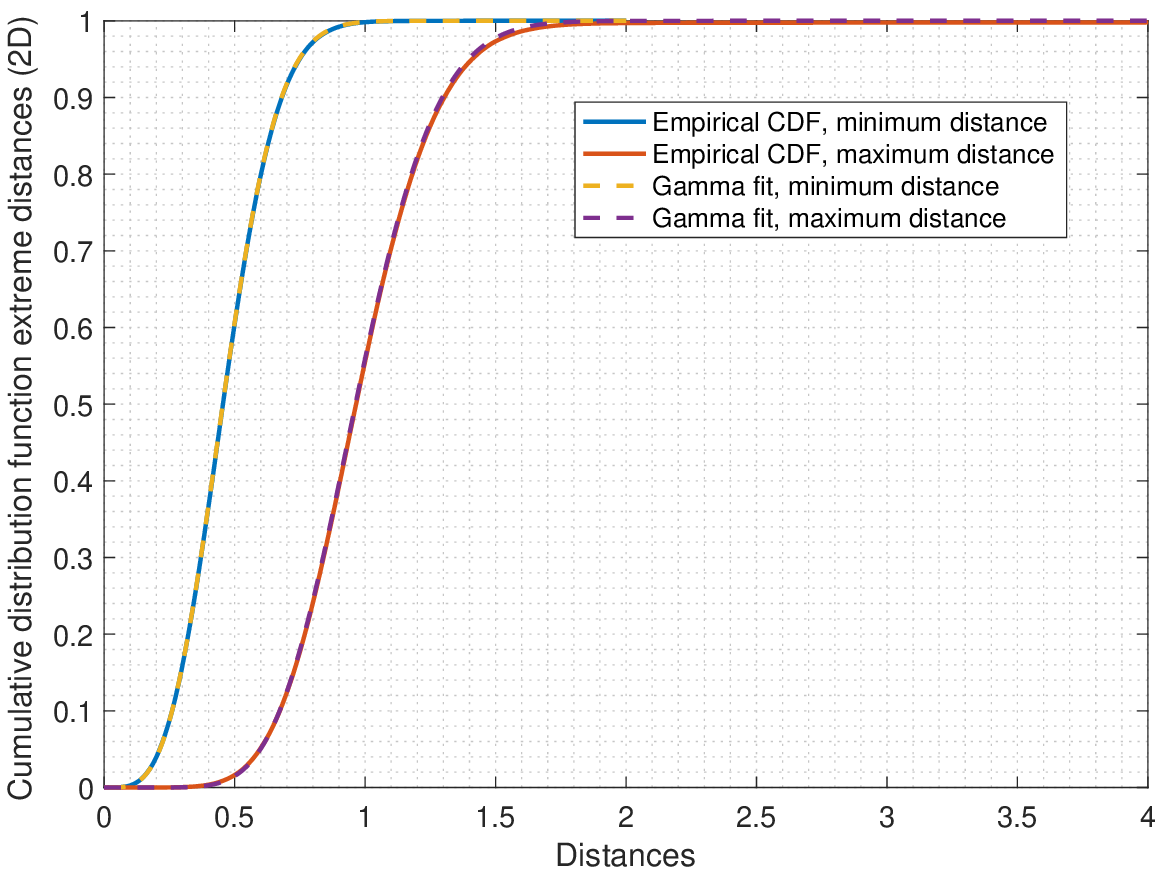}
    \caption{Cumulative distribution functions of the normalized maximum and minimum distances in 2D bounded, homogeneous Poisson Voronoi regions. Straight lines represent the empirical cumulative density function of the maximum and minimum distances that resulted from the Monte-Carlo simulation described in Algorithm \ref{alg:VoronoiTesselation_2D} with $\lambda A = 10^8$ Voronoi cells. Dashed lines represent out best-fit approximation to the Generalized Gamma distribution.} 
    \label{fig:CDFs_extremeDist_GammaFit}
\end{figure*}

Figure \ref{fig:CDFs_extremeDist_GammaFit} illustrates the results of the Monte-Carlo simulation and the best-fit approximations. The estimates of the parameters of the Generalized Gamma distribution $a$, $b$, and $c$ are presented in Table \ref{tab:GammaFitParams-GGG} with their corresponding $95\%$ confidence interval.

\begin{table}[H]
    \centering
    \caption{Maximum likelihood estimates of the Generalized Gamma distribution functions with thir 95\% confidence intervals.\\ We also append for completion the results from two other papers.}
    \begin{tabular}{|l|c|c|c|}
        \hline
        \textbf{Parameter of interest} & $\hat{\bm{a}}$ & $\hat{\bm{b}}$ & $\hat{\bm{c}}$ \\ \hline \hline
        Normalized min. Euclidean distance, $\Bar{R}_\textrm{min}$ & $ 2.176 \; (2.176, 2.176)$ & $8.446 \; (8.446, 8.446)$ & $4.005\; (4.005, 4.005)$ \\ \hline
        Normalized max. Euclidean distance, $\Bar{R}_\textrm{max}$ & $1.719\; (1.719, 1.719)$ & $5.528\; (5.527, 5.529)$ & $9.482\; (9.481, 9.482)$ \\ \hline \hline
        Normalized Euclidean distance, $\Bar{R}$ (taken from \cite{brakke1987statistics}) & 2 & $\pi$ & $4$\\ \hline
        Normalized cell Area, $\Bar{A}$ (taken from \cite{ferenc2007size}) & 1.083 & 3.011 & 2.298 \\ \hline
        \end{tabular}
    \label{tab:GammaFitParams-GGG}
\end{table}

Finally, and for completion, we derive the $m$-th order moment of the Generalized Gamma distribution as means to characterize the underlying behavior of our model, assess our model adequacy, and provide a measure for comparison to other models \cite{moller2012lectures}.
\begin{align}
    \mathbb{E}[\Bar{X}^m] & = \frac{a\,b^{c/a}}{\Gamma(c/a)}\, \int_0^\infty x^{c+m-1}\,e^{-bx^a}dx = \frac{1}{\Gamma(c/a)}\, \int_0^\infty \frac{1}{b^{m/a}}t^{(c+m)/a-1}\,e^{-t}dx = \frac{\Gamma((c+m)/a)}{b^{m/a}\,\Gamma(c/a)}.
\end{align}
From there and Table \ref{tab:GammaFitParams-GGG}, it follows that the first and second order moments of our best-fit approximations are
\begin{align*}
    \mathbb{E}[\Bar{R}_\textrm{min}] = 0.464162, \qquad & \mathbb{E}[\Bar{R}_\textrm{max}] = 0.976030, \\
    \mathbb{E}[\Bar{R}^2_\textrm{min}] = 0.241833, \qquad & \mathbb{E}[\Bar{R}^2_\textrm{max}] = 1.011695.
\end{align*}

\section{Conclusions}
\label{sec:Conclusions}
Studying distances in Voronoi tessellations reveal interesting geometric patterns, relationships, and statistical properties within the spatial arrangement of points. Our paper, focused on homogeneous Poisson point process, finds a mathematical approximation to the cumulative distribution functions governing the smallest and largest distances between the points generating the Voronoi regions and their respective vertices. We show that the best-fit approximation in terms of lowest root mean-squared error and lowest maximum absolute variation is given by the Generalized Gamma distribution. We also calculate the maximum likelihood of its parameters and derive the first and second order moments of this best-fit distribution.

As researchers continue to delve into the intricate nature of Voronoi tessellations, we anticipate that our findings will pave the way for new insights, methodologies, and applications, ultimately enriching the understanding of spatial point processes and facilitating advancements across various domains.

\bibliography{references}

\begin{thebibliography}{10}
\providecommand{\url}[1]{#1}
\csname url@samestyle\endcsname
\providecommand{\newblock}{\relax}
\providecommand{\bibinfo}[2]{#2}
\providecommand{\BIBentrySTDinterwordspacing}{\spaceskip=0pt\relax}
\providecommand{\BIBentryALTinterwordstretchfactor}{4}
\providecommand{\BIBentryALTinterwordspacing}{\spaceskip=\fontdimen2\font plus
\BIBentryALTinterwordstretchfactor\fontdimen3\font minus \fontdimen4\font\relax}
\providecommand{\BIBforeignlanguage}[2]{{%
\expandafter\ifx\csname l@#1\endcsname\relax
\typeout{** WARNING: IEEEtran.bst: No hyphenation pattern has been}%
\typeout{** loaded for the language `#1'. Using the pattern for}%
\typeout{** the default language instead.}%
\else
\language=\csname l@#1\endcsname
\fi
#2}}
\providecommand{\BIBdecl}{\relax}
\BIBdecl

\bibitem{brakke1987statistics}
K.~A. Brakke, ``Statistics of random plane {V}oronoi tessellations,'' \emph{Department of Mathematical Sciences, Susquehanna University (Manuscript 1987a)}, vol.~18, 1987.

\bibitem{muche2005poisson}
L.~Muche, ``The {P}oisson-{V}oronoi tessellation: relationships for edges,'' \emph{Advances in applied probability}, vol.~37, no.~2, pp. 279--296, 2005.

\bibitem{ferenc2007size}
J.-S. Ferenc and Z.~N{\'e}da, ``On the size distribution of {P}oisson {V}oronoi cells,'' \emph{Physica A: Statistical Mechanics and its Applications}, vol. 385, no.~2, pp. 518--526, 2007.

\bibitem{tanemura2003statistical}
M.~Tanemura, ``Statistical distributions of {P}oisson {V}oronoi cells in two and three dimensions,'' \emph{FORMA-TOKYO-}, vol.~18, no.~4, pp. 221--247, 2003.

\bibitem{edla2012clustering}
D.~R. Edla and P.~K. Jana, ``Clustering biological data using {V}oronoi diagram,'' in \emph{Advanced Computing, Networking and Security: International Conference, ADCONS 2011, Surathkal, India, December 16-18, 2011, Revised Selected Papers}.\hskip 1em plus 0.5em minus 0.4em\relax Springer, 2012, pp. 188--197.

\bibitem{stewart2010voronoi}
C.~W. Stewart and R.~van~der Ree, ``A {V}oronoi diagram based population model for social species of wildlife,'' \emph{Ecological modelling}, vol. 221, no.~12, pp. 1554--1568, 2010.

\bibitem{bherwani2021understanding}
H.~Bherwani, S.~Anjum, S.~Kumar, S.~Gautam, A.~Gupta, H.~Kumbhare, A.~Anshul, and R.~Kumar, ``Understanding {COVID}-19 transmission through {B}ayesian probabilistic modeling and {GIS}-based {V}oronoi approach: a policy perspective,'' \emph{Environment, Development and Sustainability}, vol.~23, pp. 5846--5864, 2021.

\bibitem{ramella2001finding}
M.~Ramella, W.~Boschin, D.~Fadda, and M.~Nonino, ``Finding galaxy clusters using {V}oronoi tessellations,'' \emph{Astronomy \& Astrophysics}, vol. 368, no.~3, pp. 776--786, 2001.

\bibitem{feinstein2021bionoi}
J.~Feinstein, W.~Shi, J.~Ramanujam, and M.~Brylinski, ``Bionoi: A {V}oronoi diagram-based representation of ligand-binding sites in proteins for machine learning applications,'' \emph{Protein-Ligand Interactions and Drug Design}, pp. 299--312, 2021.

\bibitem{andrews2016modeling}
J.~G. Andrews, T.~Bai, M.~N. Kulkarni, A.~Alkhateeb, A.~K. Gupta, and R.~W. Heath, ``Modeling and analyzing millimeter wave cellular systems,'' \emph{IEEE Transactions on Communications}, vol.~65, no.~1, pp. 403--430, 2016.

\bibitem{calka2014extreme}
P.~Calka and N.~Chenavier, ``Extreme values for characteristic radii of a {P}oisson-{V}oronoi tessellation,'' \emph{Extremes}, vol.~17, pp. 359--385, 2014.

\bibitem{peris2022modelling}
J.~Anguera-Peris and V.~Fodor, ``Modelling multi-cell edge video analytics,'' in \emph{IEEE International Conference on Communications}, 2022.

\bibitem{chambolle2023crystallinity}
A.~Chambolle and L.~Kreutz, ``Crystallinity of the homogenized energy density of periodic lattice systems,'' \emph{Multiscale Modeling \& Simulation}, vol.~21, no.~1, pp. 34--79, 2023.

\bibitem{jalem2018general}
R.~Jalem, M.~Nakayama, Y.~Noda, T.~Le, I.~Takeuchi, Y.~Tateyama, and H.~Yamazaki, ``A general representation scheme for crystalline solids based on voronoi-tessellation real feature values and atomic property data,'' \emph{Science and Technology of advanced MaTerialS}, vol.~19, no.~1, pp. 231--242, 2018.

\bibitem{xia2016review}
K.~Xia and G.-W. Wei, ``A review of geometric, topological and graph theory apparatuses for the modeling and analysis of biomolecular data,'' \emph{arXiv preprint arXiv:1612.01735}, 2016.

\bibitem{ma2019emergency}
Y.~Ma, W.~Xu, L.~Qin, X.~Zhao, and J.~Du, ``Emergency shelters location-allocation problem concerning uncertainty and limited resources: a multi-objective optimization with a case study in the central area of {B}eijing, {C}hina,'' \emph{Geomatics, Natural Hazards and Risk}, 2019.

\bibitem{gezer2021statistical}
F.~Gezer, R.~G. Aykroyd, and S.~Barber, ``Statistical properties of {P}oisson-{V}oronoi tessellation cells in bounded regions,'' \emph{Journal of Statistical Computation and Simulation}, vol.~91, no.~5, pp. 915--933, 2021.

\bibitem{moller2012lectures}
J.~Moller, \emph{Lectures on random {V}oronoi tessellations}.\hskip 1em plus 0.5em minus 0.4em\relax Springer Science \& Business Media, 2012, vol.~87.

\bibitem{lee1980two}
D.-T. Lee and B.~J. Schachter, ``Two algorithms for constructing a {D}elaunay triangulation,'' \emph{International Journal of Computer \& Information Sciences}, vol.~9, no.~3, pp. 219--242, 1980.

\bibitem{cai2023maximum}
Y.~Cai, ``Maximum likelihood estimates of parameters in generalized {G}amma distribution with self algorithm,'' \emph{arXiv preprint arXiv:2306.16419}, 2023.

\end{thebibliography}
\bibliographystyle{IEEEtran}

\appendix

\begin{algorithm}
\label{alg:VoronoiTesselation_2D}
\DontPrintSemicolon
\caption{Maximum and minimum distances in interior cells of Poisson Voronoi tessellations (2D)}
\KwIn{Area of the planar space $A$ and density of seed points per unit area $\lambda$}
\KwOut{Cumulative distribution function of the maximum and minimum distances of a sample realization of an homogeneous Poisson Point process}

Generate a Poisson random variable $ N \leftarrow \text{Poisson}(\lambda L)$, and define three sets $ \mathcal{P} \leftarrow \text{Matrix}(\textit{size} : N \times 2)$, $ R_\mathrm{min} \leftarrow \text{Vector}(\textit{size} : N)$, $ R_\mathrm{max} \leftarrow \text{Vector}(\textit{size} : N)$, and a dynamic set $ \ell \leftarrow \text{Vector}^*(\textit{size} : N)$\;

Generate $N$ seed points $p_i$ uniformly at random on the space $[0,\sqrt{A}]^2$ and store them in $\mathcal{P}$\;

Employ the Delanuay triangulation \cite[Section~3.2]{lee1980two} to obtain all Voronoi regions $\mathcal{V}(p_i),\;\forall p_i\in\mathcal{P}$ and all edges \;

Construct a Voronoi diagram from the results of the Delanuay triangulation and obtain all vertices\;

\For{$i\leftarrow 1$ \KwTo $N$}{
    Find all vertices associated to $\mathcal{V}(p_i)$ and store the result in $\ell$\;
}

\For{$i\leftarrow 1$ \KwTo $N$}{
    Calculate all the Euclidean distances between $p_i$ and all its associated vertices $x\in\ell(p_i)$, $r_i =  \| x - p_i\|_2$\;
    Store the maximum and minimum distances in $R_\mathrm{min}$ and $R_\mathrm{max}$, respectively\;
}
Calculate the normalized maximum and minimum distances $\Bar{R}_\mathrm{min} = \sqrt{\lambda} R_\mathrm{min}$ and $\Bar{R}_\mathrm{max} = \sqrt{\lambda} R_\mathrm{max}$\;
Plot the empirical CDF of $\Bar{R}_\mathrm{min}$ and $\Bar{R}_\mathrm{max}$ and their best-fit approximations\;
Calculate the max. likelihood and confidence intervals of the parameters of the Generalized Gamma distribution \cite{cai2023maximum}\;

\end{algorithm}

\begin{table}[H]
    \centering
    \caption{Fit approximations to the empirical normalized maximum and minimum distances.\\The resulting fits are compared in terms of the root mean-squared error and the maximum absolute variation, defined as the squared root of the 2-norm and the $\infty$-norm of the difference between the empirical and the fit distributions, respectively.}
    \begin{tabular}{|c|c|c||c|c|c|}
        \hline
        \multicolumn{3}{|c||}{\textbf{Minimum normalized Euclidean distance}} & \multicolumn{3}{c|}{\textbf{Maximum normalized Euclidean distance}} \\ \hline
        \textit{Distributions} & \textit{RMSE} & \textit{Max. abs. variation} & \textit{Distributions} & \textit{RMSE} & \textit{Max. abs. variation}\\ \hline \hline
        Generalized Gamma & $0.0001345$ & $0.00205$ & Generalized Gamma & $0.0006213$ & $0.00487$\\ \hline
        Gamma & $0.0094499$ & $0.01700$ & Gamma & $0.0046578$ & $0.00876$ \\ \hline
        Log-normal & $0.0179480$ & $0.03240$ & Log-normal & $0.0105930$ & $0.01902$ \\ \hline
        Rayleigh & $0.0759310$ & $0.90862$ & Rayleigh & $0.1389700$ & $0.05735$ \\ \hline
        Weibull  & $0.0066456$ & $0.01218$ & Weibull  & $0.0140030$ & $0.02822$ \\ \hline
        \end{tabular}
    \label{tab:GammaFitParams}
\end{table}

\end{document}